\documentclass[a4paper,oneside,10pt]{article}%
\usepackage{amsmath}
\usepackage{amsfonts}
\usepackage{amssymb}
\usepackage{appendix}
\usepackage{graphicx}
\usepackage[square,numbers,sort&compress]{natbib}
\usepackage{color}%
\providecommand{\U}[1]{\protect\rule{.1in}{.1in}}

\newtheorem{theorem}{Theorem}[section]

\newtheorem{definition}[theorem]{Definition}
\newtheorem{assumption}{Assumption}

\newtheorem{lemma}[theorem]{Lemma}

\newtheorem{remark}[theorem]{Remark}

\numberwithin{equation}{section}

\newcommand{\E}{{\mathbb E}}
\newcommand{\R}{{\mathbb R}}
\newcommand{\pf}{\noindent\textbf{Proof:} }
\newcommand{\eof}{\hfill{$\Box$}}

\newcommand{\cM}{\ensuremath{\mathcal{M}}}

\usepackage{hyperref}
\hypersetup{
colorlinks=true, 
breaklinks=true, 
urlcolor= blue, 
linkcolor= red, 
citecolor=blue, 
pdftitle={}, 
pdfauthor={}, 
pdfsubject={} 
}

\begin{document}

\title{Stochastic linear-quadratic control with a jump and regime switching on a random horizon}
\author{Ying Hu \thanks{Univ Rennes, CNRS, IRMAR-UMR 6625, F-35000 Rennes, France. Partially supported by Lebesgue
Center of Mathematics \textquotedblleft Investissements d'avenir\textquotedblright program-ANR-11-LABX-0020-01, ANR CAESARS
(No.~15-CE05-0024) and ANR MFG (No.~16-CE40-0015-01). Email:
\url{ying.hu@univ-rennes1.fr }}
\and Xiaomin Shi\thanks{School of Mathematics and
Quantitative Economics, Shandong University of Finance and Economics, Jinan
250100, China. Partially supported by NSFC (No.~11801315), NSF of Shandong Province (No.~ZR2018QA001). Email: \url{shixm@mail.sdu.edu.cn}}
\and Zuo Quan Xu\thanks{Department of Applied Mathematics, The Hong Kong Polytechnic University, Kowloon, Hong Kong.
Partially supported by NSFC (No.~11971409), Hong Kong RGC
(GRF No.~15202421), The PolyU-SDU Joint Research Center on Financial Mathematics and the CAS AMSS-PolyU Joint Laboratory of Applied Mathematics, The Hong Kong Polytechnic University. Email: \url{maxu@polyu.edu.hk}}}
\maketitle

\textbf{Abstract}.
In this paper, we study a stochastic linear-quadratic control problem with random coefficients and regime switching on a horizon $[0,T\wedge\tau]$, where $\tau$ is a given random jump time for the underlying state process and $T$ is a constant.
We obtain an explicit optimal state feedback control and explicit optimal cost value
by solving a system of stochastic Riccati equations (SREs) with jumps on $[0,T\wedge\tau]$. By the decomposition approach stemming from filtration enlargement theory, we express the solution of the system of SREs with jumps in terms of another system of SREs involving only Brownian filtration on the deterministic horizon $[0,T]$. Solving the latter system is the key theoretical contribution of this paper and
we establish this for three different cases, one of which seems to be new in the literature. These results are then applied to study a mean-variance hedging problem with random parameters that depend on both Brownian motion and Markov chain. The optimal portfolio and optimal value are presented in closed forms with the aid of a system of linear backward stochastic differential equations with jumps and unbounded coefficients in addition to the SREs with jumps.



\addcontentsline{toc}{section}{\hspace*{1.8em}Abstract}

\section{Introduction}
The stochastic linear-quadratic (LQ, for short) optimal control has known important developments since the pioneering works of Wonham \cite{Wo} and Bismut \cite{Bi} who studied stochastic LQ problem with deterministic and random coefficients respectively. Kohlmann and Zhou \cite{KZ} established the relationship between stochastic LQ problems and backward stochastic differential equations (BSDEs). Chen, Li and Zhou \cite{CLZ} studied the indefinite stochastic LQ problem which is different significantly from its deterministic counterpart. Meanwhile, it has wide applications in many fields, such as risk management, optimal investment or mean-variance portfolio selection; see, e.g., \cite{HZ, KT, LZL, Lim, LZ, ZL, ZY}.

All the aforementioned literature focused on deterministic time horizon. However the controller in reality may only access the model before the realization of some random event, such as before the default time in credit risk theory and death time in actuarial science. Therefore, it is of great theoretical importance to study control problems with random time of eventual exit.

Yu \cite{Yu} and Lv, Wu and Yu \cite{LY} studied continuous time mean-variance portfolio selection problems with random horizon in complete and incomplete markets, respectively, where the state process is continuous.
They assumed the conditional distribution of the eventual exit follows an It\^{o} process driven by the same Brownian motion as in risky assets dynamics.
In practice, there is another situation that is the time of eventual exit $\tau$ arrives by surprise, i.e. $\tau$ is a totally inaccessible random time for the reference filtration, meanwhile the state process may have a jump at the time $\tau$.
The theory of enlargement of filtration happens to be a powerful tool for modeling such random time. Please refer to Bielecki and Rutkowski \cite{BR} for a systematic account on this subject. Along this line, Pham \cite{Pham} investigated a general stochastic control problem under a progressive enlargement of filtration and proved a decomposition of the
original stochastic control problem under the global filtration into classical stochastic control problems under the reference filtration. Jeanblanc, Mastrolia, Possamai and Reveillac \cite{JMPR} studied an exponential utility maximization problem with bounded random horizon. 
Neither \cite{Pham} nor \cite{JMPR} considered jumps in the state (or wealth) processes. Kharroubi, Lim and Ngoupeyou \cite{KLN} studied a mean-variance hedging problem with jumps on random horizon. The main approach used in \cite{JMPR} and \cite{KLN} is to reduce their problems to the analysis of the solvability of some BSDEs with jumps.

We generalize Kharroubi, Lim and Ngoupeyou's \cite{KLN} model to a general stochastic linear-quadratic control problem with regime switching.
All the coefficients and the control variable (portfolio) in \cite{KLN} are one dimensional. By contrast,
both the control variable and the Brownian motion in our model are multi-dimensional. Moreover, we will solve the stochastic LQ problems and the associate stochastic Riccati equations (SREs) with jumps
for three different cases: one standard case and two singular cases with one of which seems to be new in the literature, while only one singular case was considered in \cite{KLN}.

On the other hand,
a Markov chain is usually adopted to reflect the model status or the random market environment. For instance, Zhou and Yin \cite{ZY} studied a mean-variance portfolio selection problem, where the coefficients are assumed to be deterministic functions of time for each given regime. To better reflect randomness of the market environment, all the coefficients in our model are assumed to be stochastic processes for each given regime.

As is well known, the SREs plays a crucial role in representing the optimal controls for stochastic LQ control problems. It is Tang \cite{Ta} who established the existence and uniqueness results for SREs with uniformly definite coefficients. While there were only partial results for the solvability of indefinite SREs so far; see, e.g., \cite{HZ03}. Zhang, Dong and Meng \cite{ZDM} made a great progress in solving stochastic LQ control and related SREs with jumps with uniformly definite control weight by inverse flow technique.

The main theoretical contribution of this paper is the solvability of the associated systems of SREs with jumps.
Due to the presence of regime switching, the SREs with jumps in our model is actually a system of BSDEs coupled through the generator of the Markov chain.
Different from the inverse flow method used in \cite{ZDM}, we will establish the solvability of our SREs with jumps from the point of view of BSDE theory directly. Inspired by \cite{KLN}, we construct the solution to the SREs on random horizon with jumps from another SREs on deterministic finite horizon without jumps. Using the model of {our previous work} \cite{HSX}, we obtain the solvability of the latter SREs. Different with \cite{HSX}, new terms emerge in the SREs in this paper so that we have to establish nonnegative or uniformly positive lower bounds for their solution.

The rest part of this paper is organized as follows. Section 2 presents the general framework and assumptions of stochastic LQ problem on random horizon. In Section 3, we establish the solvability of the system of SREs with jumps. The optimal feedback control and optimal value are presented in Section 4 with the aid of SREs with jumps. In Section 5, the general results are applied to solve a mean-variance hedging problem.

\section{Problem formulation}
Let $(\Omega, \mathcal H, \mathbb{P})$ be a fixed complete probability space on which is defined a standard $n$-dimensional Brownian motion {$\{W_t\}_{t\geq 0}$.
Define $\mathbb{F}=(\mathcal{F}_t)_{t\geq0}$, $\mathcal F_t=\sigma\{W_s: 0\leq s\leq t\}\bigvee\mathcal{N}$, where $\mathcal{N}$ is the totality of all the $\mathbb{P}$-null sets of $\mathcal{H}$.
In this probability space, a random time $\tau$ is given, which represents, for example, the default time in credit or counterparty risk models, or a death time in actuarial models. The random time $\tau$ is not assumed to be an $\mathbb{F}$-stopping time. We therefore use in the sequel the standard approach of progressive enlargement of filtrations. Let $\mathbb{G}$ be the smallest right continuous extension of $\mathbb{F}$ that turns $\tau$ into a $\mathbb{G}$-stopping time. More precisely, $\mathbb{G}:=(\mathcal{G}_t)_{t\geq 0}$ is defined by
\[
\mathcal{G}_t:=\bigcap_{\varepsilon>0}\widetilde{\mathcal{G}}_{t+\varepsilon}, \ \mbox{ for all $t\geq0$},
\]
where $\widetilde{\mathcal{G}}_{s}:=\mathcal{F}_s\vee\sigma(\mathbf{1}_{\tau\leq u}; 0\leq u\leq s)$, for all $s\geq 0$.

On $(\Omega, \mathcal H, \mathbb{P})$, there is a continuous-time stationary Markov chain $\{\alpha_t\}_{t\geq 0}$ valued in a finite state space $\mathcal M=\{1,2,...,\ell\}$ such that $\{\alpha_t\}_{t\geq 0}$ is independent of $\{W_t\}_{t\geq 0}$ and $\tau$. The Markov chain $\{\alpha_t\}_{t\geq 0}$ has a generator $Q=(q_{ij})_{\ell\times \ell}$ with $q_{ij}\geq0$ for $i\neq j$ and $\sum_{j=1}^\ell q_{ij}=0$ for every $i,j\in\cM$.
Denote by $\mathbb{F}^{\alpha}=(\mathcal{F}_t^{\alpha})_{t\geq0}$ the natural filtration generated by $\alpha$,
and
$\mathbb{H}:=(\mathcal{H}_t)_{t\geq 0}$, where
$\mathcal{H}_t:=\mathcal{G}_t\vee\mathcal{F}_t^{\alpha}$.

We denote by $\mathcal{P}(\mathbb{F})$ the $\sigma$-algebra of $\mathbb{F}$-predictable subsets
of $\Omega\times\mathbb{R}^{+}$, i.e. the $\sigma$-algebra generated by the left-continuous $\mathbb{F}$-adapted
processes. Define $\mathcal{P}(\mathbb{G})$ and $\mathcal{P}(\mathbb{H})$ in a similar way.

We introduce the following notation:
\[%
\begin{array}
[c]{l}%
L^{\infty}_{\mathcal{F}_T}(\Omega)=\Big\{\xi:\Omega\rightarrow
\mathbb{R}\;\Big|\;\xi\mbox { is }\mathcal{F}_{T}\mbox{-measurable, and essentially bounded}\Big\},\\
L^{2}_{\mathbb F}(0,T;\mathbb{R})=\Big\{\phi:[0,T]\times\Omega\rightarrow
\mathbb{R}\;\Big|\;(\phi_{t})_{0\leq t\leq T}\mbox{ is }\mathcal{P}(\mathbb{F})\mbox{-measurable such that }\E\int_{0}^{T}|\phi_{t}|^{2}dt<\infty
\Big\},\\
L^{\infty}_{\mathbb{F}}(0,T;\mathbb{R})=\Big\{\phi:[0,T]\times\Omega
\rightarrow\mathbb{R}\;\Big|\;(\phi_{t})_{0\leq t\leq T}\mbox{ is }\mathcal{P}(\mathbb{F})\mbox{-measurable essentially bounded} \Big\},\\
S^{\infty}_{\mathbb{F}}(0,T;\mathbb{R})=\Big\{\phi:[0,T]\times\Omega
\rightarrow\mathbb{R}\;\Big|\;(\phi_{t})_{0\leq t\leq T}\mbox{ is } \ \mbox{c}\grave{\mathrm{a}}\mbox{d-l}\grave{\mathrm{a}}\mbox{g} \ \mathbb{F}\mbox{-adapted essentially bounded}\Big\}.
\end{array}
\]
These definitions are generalized in the obvious way to the cases that $\mathbb{F}$ is replaced by $\mathbb{G}$, $\mathbb{H}$, $[0,T]$ by any random time $[\tau,\iota]$ and $\mathbb{R}$ by $\mathbb{R}^n$, $\mathbb{R}^{n\times m}$ or $\mathbb{S}^n$, where $\mathbb{S}^n$ is the set of symmetric $n\times n$ real matrices. If $M\in\mathbb{S}^n$ is positive definite (resp. positive semidefinite), we write $M>0$ (resp. $M\geq 0$).
In our argument, $t$, $\omega$, ``almost surely'' and ``almost everywhere'', may be suppressed for simplicity in many circumstances if no confusion occurs. In this paper the integral $\int_s^t$ stands for $\int_{(s,t]}$.

The following assumption is classical in the filtration enlargement theory.

\begin{assumption}\label{assu-2}
Any $\mathbb{F}$-martingale is a $\mathbb{G}$-martingale.
\end{assumption}

\begin{remark}
Assumption \ref{assu-2} is known as immersion between $\mathbb{F}$ and $\mathbb{G}$ in the filtration enlargement theory. Under Assumption \ref{assu-2}, the $\sigma$-field $\mathcal{F}_\infty$ and $\mathcal{G}_t$ are conditionally independent given $\mathcal{F}_t$, and the process {$\{W_t\}_{t\geq 0}$} remains a $\mathbb{G}$-Brownian motion. Please see Theorem 3.2 in Aksamit and Jeanblanc \cite{AJ} and Remark 4.1 in Kharroubi and Lim \cite{KL}. Since $\mathbb{F}^{\alpha}$ is independent of $\mathbb{G}$, we further conclude that $\{W_t\}_{t\geq 0}$ remains an $\mathbb{H}$-Brownian motion.
\end{remark}

Since $\mathcal{F}_t^{\alpha}$ is independent of $\mathcal{G}_t$, under Assumption \ref{assu-2}, the stochastic integral $\int_0^tX_sdW_s$ is well defined for all $\mathcal{P}(\mathbb{H})$-measurable process $X$ such that $\int_0^t|X_s|^2ds<\infty$.

\begin{assumption}\label{assu-1}
The process $N_{\cdot}:=\mathbf{1}_{\tau\leq \cdot}$ admits an $\mathbb{F}$-compensator of the form $\int_0^{\cdot\wedge\tau}\lambda_sds$, i.e. $$M_{\cdot}:=N_{\cdot}-\int_0^{\cdot\wedge\tau}\lambda_sds$$ is a $\mathbb{G}$-martingale, where $\lambda$ is a bounded nonnegative $\mathcal{P}(\mathbb{F})$-measurable process.
\end{assumption}

We remark that
\begin{align*}
M_t
&=N_t-\int_0^{t}\lambda_s\mathbf{1}_{s\leq\tau}ds
=N_t-\int_0^{t}\lambda_s^{\mathbb{G}}ds, \quad t\geq0,
\end{align*}
where $$\lambda_s^{\mathbb{G}}:=\lambda_s\mathbf{1}_{s\leq\tau}=\lambda_s(1-N_{s-})\geq 0, \quad s\geq0,$$
is a $\mathcal{P}(\mathbb{G})$-measurable process.

We now introduce the following scalar-valued linear stochastic differential equation (SDE):
\begin{align}
\label{state}
\begin{cases}
dX_t=\left[A^{\alpha_t}_t X_{t-}+u_t'B^{\alpha_t}_t\right]dt
+\left[C^{\alpha_t}_t X_{t-}+D^{\alpha_t}_tu_t\right]'dW_t
+[E^{\alpha_t}_tX_{t-}+u_t'F^{\alpha_t}_t]dM_t, \ t\in[0, T\wedge\tau], \\
X_0=x, \ \alpha_0=i_0,
\end{cases}
\end{align}
where $T$ is a positive constant, $A^i, \ B^i, \ C^i, \ D^i, \ E^i, \ F^i$ are all $\mathcal{P}(\mathbb{F})$-measurable processes of suitable sizes for all $i\in\cM$, $x\in\mathbb{R}$ and $i_0\in\cM$ are known.

The class of admissible controls is defined as the set
\begin{align*}
\mathcal{U}:= L^{2}_{\mathbb H}(0,T\wedge\tau;\mathbb{R}^m).
\end{align*}
If $u(\cdot)\in\mathcal{U}$ and $X(\cdot)$ is the associated solution of \eqref{state}, then we refer to $(X(\cdot), u(\cdot))$ as an admissible pair.

Let us now state our stochastic LQ problem as follows:
\begin{align}
\begin{cases}
\mathrm{Minimize} &\ J(x, i_0, u(\cdot))\\
\mbox{subject to} &\ (X, u) \mbox{ is an admissible pair for} \ \eqref{state},
\end{cases}
\label{LQ}%
\end{align}
where the cost functional is given as the following quadratic form
\begin{align}\label{costfunctional}
J(x, i_0, u(\cdot)):=\mathbb{E}\left\{\int_0^{T\wedge\tau}\Big(Q_t^{\alpha_t}X_t^2+u_t'R_t^{\alpha_t}u_t\Big)dt
+G^{\alpha_{T\wedge\tau}}X_{T\wedge\tau}^2\right\}.
\end{align}

For $(x, i_0)\in\R\times\cM$, the problem \eqref{LQ} is said to be well-posed, if
\begin{align*}
V(x,i_0):=\inf_{u\in\mathcal{U}}J(x, i_0, u)>-\infty,
\end{align*}
where $V(x,i_0)$ is called its optimal value; and
the problem is said to be solvable, if there exists a control $u^*\in\mathcal{U}$ such that
\begin{align*}
J(x, i_0, u^*) =V(x,i_0)>-\infty,
\end{align*}
in which case, $u^*$ is called an optimal control for the problem \eqref{LQ}.

We introduce a decomposition result for $\mathcal{P}(\mathbb{G})$-measurable processes. Please refer to Proposition 2.1 in \cite{KLN}.
\begin{lemma}
\label{gtof}
Any $\mathcal{P}(\mathbb{G})$-measurable process $Y=(Y_t)_{t\geq0}$ can be represented as
\[
Y_t=Y_t^b\mathbf{1}_{t\leq\tau}+Y_t^a(\tau)\mathbf{1}_{t>\tau},\quad t\geq 0,
\]
where $Y^b$ is $\mathcal{P}(\mathbb{F})$-measurable and $Y^a$ is $\mathcal{P}(\mathbb{F})\otimes\mathcal{B}(\mathbb{R}^{+})$-measurable. 
\end{lemma}

\begin{remark}
By Lemma \ref{gtof}, it is reasonable to assume that the coefficients $A^i,\ B^i, \ C^i, \ D^i, \ E^i, \ F^i, \ Q^i, \ R^i$ are $\mathcal{P}(\mathbb{F})$-measurable, because the LQ problem \eqref{LQ} depends only on their values in the interval $[0,T\wedge\tau]$.
\end{remark}

\begin{assumption} \label{assu1}
For all $i\in\cM$,
\begin{align*}
\begin{cases}
A^i\in L_{\mathbb{F}}^\infty(0, T;\mathbb{R}), \
B^i\in L_{\mathbb{F}}^\infty(0, T;\mathbb{R}^m), \
C^i \in L_{\mathbb{F}}^\infty(0, T;\mathbb{R}^n), \\
D^i\in L_{\mathbb{F}}^\infty(0, T;\mathbb{R}^{n\times m}), \
E^i\in L_{\mathbb{F}}^\infty(0,T;\mathbb{R}),\
F^i\in L_{\mathbb{F}}^\infty(0,T;\mathbb{R}^m),\\
Q^i\in L_{\mathbb{F}}^\infty(0, T;\mathbb{R}), \
R^i\in L_{\mathbb{F}}^\infty(0, T;\mathbb{S}^m),
\end{cases}
\end{align*}
and
$G^i$ is a bounded $\mathcal{G}_{T\wedge\tau}$-measurable random variable of the form
\begin{align*}
G^i=G^{b,i}\mathbf{1}_{T<\tau}+G^{a,i}_{\tau}\mathbf{1}_{T\geq\tau},
\end{align*}
where $G^{b,i}\in L^{\infty}_{\mathcal{F}_T}(\Omega)$ and
$G^{a,i}\in L^\infty_{\mathbb{F}}(0,T;\mathbb{R})$.
\end{assumption}
Here the superscript ``$b$" (resp. ``$a$") stands for ``before the jump" (resp. ``after the jump").

Throughout the paper, the above Assumptions \ref{assu-2}, \ref{assu-1} and \ref{assu1} are always, implicitly or explicitly, put in force.

\begin{remark}
There are three kinds of uncertainties in the problem \eqref{LQ}. One stems from the randomness of the Brownian motion, the second one comes from the random time, and the last one comes from the Markov chain.
\end{remark}

This paper assumes the coefficients in the cost functional \eqref{costfunctional} fulfill at least one of the following conditions, so that the problem \eqref{LQ} is well-posed with a nonnegative optimal value.

\begin{assumption} [Standard Case]
\label{assu3}
$Q^i\geq0$, $G^i\geq0$ and there exists a constant $\delta>0$ such that $R^i\geq \delta I_{m}$, for all $i\in\cM$, where $I_{m}$ denotes the $m$-dimensional identity matrix.
\end{assumption}

\begin{assumption}[Singular Case I]
\label{assu2}
$Q^i\geq0$, $R^i\geq0$, and
there exists a constant $\delta>0$ such that $G^{i}\geq\delta$ and $(D^i)'D^i\geq\delta I_{m}$, for all $i\in\cM$.
\end{assumption}

\begin{assumption}[Singular Case II]
\label{assu4}
$m=1$, $Q^i\geq0$, $R^i\geq0$, $G^{b,i}\geq0$ and
there exists a constant $\delta>0$ such that $G^{a,i}\geq\delta$ and {$\lambda(F^i)^2\geq\delta$}, for all $i\in\cM$.
\end{assumption}
Here ``singular" means that the control weight matrix $R^i$ in the cost functional \eqref{costfunctional} could be possibly a singular matrix.

\section{Solvability of stochastic Riccati equations}
In order to solve the problem \eqref{LQ}, we introduce the following system of ($\ell$-dimensional) BSDEs with jumps:
\begin{align}
\begin{cases}
P^i_t=G^i+\int_{t\wedge\tau}^{T\wedge\tau}\Big[(2A^i_s+(C^i_s)'C^i_s)P^i_{s-}
+\lambda^{\mathbb{G}}_s(E^i_s)^2(P^i_{s-}+U^i_s)
+2(C^i_s)'\Lambda^i_s+2\lambda^{\mathbb{G}}_sE^i_sU^i_s+Q^i_s+\sum\limits_{j=1}^\ell q_{ij}P^j_{s-}\\
\qquad\quad-\mathcal{N}(s,P^i_{s-},\Lambda^i_s,U^i_s,i)'\mathcal{R}(s,P^i_{s-},U^i_s,i)^{-1}
\mathcal{N}(s, P^i_{s-},\Lambda^i_s,U^i_s,i)\Big]ds
-\int_{t\wedge\tau}^{T\wedge\tau}(\Lambda^i_s)'dW_s-\int_{t\wedge\tau}^{T\wedge\tau}U^i_sdM_s,\\
\mathcal{R}(s,P^i_{s-},U^i, i)>0, \ \mbox{for $s\leq T\wedge\tau$ and all $i\in\cM$},
\end{cases}
\label{P}
\end{align}
where for any $(s,P,\Lambda,U)\in[0,T\wedge\tau]\times\mathbb{R}\times\mathbb{R}^n\times\mathbb{R}$,
\begin{align}\label{R}
\mathcal{N}(s,P,\Lambda,U,i)&:=P((D^i_s)' C^i_s+B^i_s)+(D^i_s)'\Lambda
+\lambda^{\mathbb{G}}_sE^i_sF^i_s(P+U)+\lambda^{\mathbb{G}}_sF^i_sU,\nonumber\\
\mathcal{R}(s,P,U,i) &:=R^i_s+P(D^i_s)'D^i_s+\lambda^{\mathbb{G}}_s(P+U)F^i_s(F^i_s)'.
\end{align}
The system of BSDEs \eqref{P} is referred as system of stochastic Riccati equations with jumps. Please note that the $\ell$ equations in \eqref{P} are coupled through $\sum_{j=1}^\ell q_{ij}P^j_{s-}$.

\begin{definition}
A vector process $(P^i, \Lambda^i, U^i)_{i\in\cM}$ is called a solution of the system of BSDEs \eqref{P}, if it satisfies \eqref{P}, and $(P^i, \Lambda^i, U^i)\in S^\infty_{\mathbb{G}}(0, T; \mathbb {R})\times L^{2}_{\mathbb{G}}(0, T;\mathbb{R}^{n})\times L^{\infty}_{\mathbb{G}}(0,T;\mathbb{R})$ for all $i\in\cM$.
Furthermore, it
is called nonnegative if $P^i\geq0, \ t\in[0,T]$, a.s.; called uniformly positive if $P^i\geq c, \ t\in[0,T]$, a.s.; and called positive if $P^i\geq0, \ P^i+U^i\geq c, \ t\in[0,T]$, a.s. where $c>0$ is some deterministic constant.
\end{definition}

We will construct a solution of \eqref{P} through another BSDE driven only by $W$, brought the idea from \cite{KLN}. We briefly recall the idea for the reader's convenience.

For any $\mathcal{G}_{T\wedge\tau}$-measurable random variables $\xi^i$ of the form
\begin{align*}
\xi^i=\xi^{b,i}\mathbf{1}_{T<\tau}+\xi^{a,i}_{\tau}\mathbf{1}_{T\geq\tau},
\end{align*}
where $\xi^{b,i}\in L^{\infty}_{\mathcal{F}_T}(\Omega)$ and $\xi^{a,i}\in L^{\infty}_{\mathbb{F}}(0,T;\mathbb{R})$,
and any $\mathcal{P}(\mathbb{G})\otimes\mathcal{B}(\mathbb{R}^\ell)\otimes\mathcal{B}(\mathbb{R})\otimes\mathcal{B}(\mathbb{R})$-measurable function $F:\Omega\times[0,T]\times\mathbb{R}^{\ell}\times\mathbb{R}\times\mathbb{R}\rightarrow\mathbb{R}$, consider the following BSDE driven by $W$ and $N$ (given in Assumption \ref{assu-1}) on horizon $[0,T\wedge\tau]$:
\begin{equation}
P^{i}_t=\xi^{i}+\int_{t\wedge\tau}^{T\wedge\tau} F_s(P_{s-},\Lambda^{i}_s,U^i_s)ds
-\int_{t\wedge\tau}^{T\wedge\tau}(\Lambda^{i}_s)'dW_s
-\int_{t\wedge\tau}^{T\wedge\tau}U_s^idN_s, \ \mbox{for all} \ i\in\cM,
\label{genejump}
\end{equation}
where $P=(P^{1}, ..., P^{\ell})'$.
Similar to \cite{KLN}, we will construct a solution to \eqref{genejump} from BSDEs without jumps in the Brownian filtration.

By Lemma \ref{gtof}, there exists a $\mathcal{P}(\mathbb{F})\otimes\mathcal{B}(\mathbb{R}^\ell)\otimes\mathcal{B}(\mathbb{R})\otimes\mathcal{B}(\mathbb{R})$-measurable function $F^b$ such that
\begin{align*}
F_t(\cdot,\cdot,\cdot)\mathbf{1}_{t\leq\tau}=F^b_t(\cdot,\cdot,\cdot)\mathbf{1}_{t\leq\tau}, \ t\geq0.
\end{align*}
The following lemma is a generalization of Theorem 4.1 in \cite{KLN} to systems of BSDEs, but the proof follows the same lines as in \cite{KLN}, so we omit it.
\begin{lemma}\label{link}
Assume that $(P^{b,i},\Lambda^{b,i})_{i\in\cM}$ is a solution to the following BSDE:
\begin{align*}
P^{b,i}_t=\xi^{b,i}+\int_t^T F^b_s(P^{b}_s,\Lambda^{b,i}_s,G^{a,i}_s-P^{b,i}_s)ds-\int_t^T(\Lambda^{b,i}_s)'dW_s, \ i\in\cM,
\end{align*}
where $P^b=(P^{b,1}, ..., P^{b,\ell})'$. Then BSDE \eqref{genejump} admits a solution
$(P^{i},\Lambda^{i}, U^i)_{i\in\cM}$ given by
\begin{align*}
&P^i_t=P^{b,i}_t\mathbf{1}_{t<\tau}+\xi^{a,i}_\tau\mathbf{1}_{t\geq\tau},\\
&\Lambda^i_t=\Lambda^{b,i}_t\mathbf{1}_{t\leq\tau},\\
&U^i_t=(\xi^{a,i}_t-P^{b,i}_t)\mathbf{1}_{t\leq\tau}, \ t\in[0,T], \ i\in\cM.
\end{align*}
\end{lemma}

\bigskip

Inspired by Lemma \ref{link}, we introduce the following system of BSDEs \emph{without jumps} on the deterministic horizon $[0,T]$:
\begin{align}
\label{PBM}
\begin{cases}
P^{b,i}_t=G^{b,i}+\int_{t}^{T}\Big[(2A^i_s+(C^i_s)'C^i_s)P^{b,i}_s+\lambda_s(E^i_s)^2G^{a,i}_s
+2(C^i_s)'\Lambda^{b,i}_s\\
\qquad\qquad\qquad+2\lambda_s E^i_s(G^{a,i}_s-P^{b,i}_s)+Q^i_s+\lambda_s(G^{a,i}_s-P^{b,i}_s)+\sum\limits_{j=1}^\ell q_{ij}P^{b,j}_s\\
\qquad\qquad\qquad-\widehat{\mathcal{N}}(s,P^{b,i}_s,\Lambda^{b,i}_s,i)'\widehat{\mathcal{R}}(s,P^{b,i}_s,i)^{-1}
\widehat{\mathcal{N}}(s,P^{b,i}_s,\Lambda^{b,i}_s,i)\Big]ds
-\int_{t}^{T}(\Lambda^{b,i}_s)'dW_s,\\
\qquad=G^{b,i}+\int_{t}^{T}\Big[(2A^i_s+(C^i_s)'C^i_s-2\lambda_s E^i_s-\lambda_s)P^{b,i}_s
+2(C^i_s)'\Lambda^{b,i}_s+\lambda_s G^{a,i}_s(E^i_s+1)^2 +Q^i_s\\
\qquad\qquad+\sum\limits_{j=1}^\ell q_{ij}P^{b,j}_s-\widehat{\mathcal{N}}(s,P^{b,i}_s,\Lambda^{b,i}_s,i)'\widehat{\mathcal{R}}(s,P^{b,i}_s,i)^{-1}
\widehat{\mathcal{N}}(s,P^{b,i},\Lambda^{b,i},i)\Big]ds
-\int_{t}^{T}(\Lambda^{b,i}_s)'dW_s,\\
\widehat{\mathcal{R}}(s,P^{b,i}_s,i)>0,\ \mbox{for $s\leq T$ and all $i\in\cM$},
\end{cases}
\end{align}
where
\begin{align*}
\widehat{\mathcal{N}}(s,P,\Lambda,i):&=P((D^i_s)' C^i_s+B^i_s)+(D^i_s)'\Lambda
+\lambda_s G^{a,i}_sE^i_sF^i_s+\lambda_s(G^{a,i}_s-P)F^i_s\nonumber\\
&=P((D^i_s)' C^i_s+B^i_s-\lambda_s F^i_s)+(D^i_s)'\Lambda+\lambda_s G^{a,i}_s(E^i_s+1)F^i_s,\nonumber\\
\widehat{\mathcal{R}}(s,P,i):&=R^i_s+P(D^i_s)'D^i_s+\lambda_s G^{a,i}_sF^i_s(F^i_s)'.
\end{align*}

\begin{definition}
A vector process $(P^i, \Lambda^i)_{i\in\cM}$ is called a solution of the system of BSDEs \eqref{PBM}, if it satisfies \eqref{PBM}, and $(P^i, \Lambda^i)\in L^\infty_{\mathbb{F}}(0, T; \mathbb {R})\times L^{2}_{\mathbb F}(0, T;\mathbb{R}^{n})$ for all $i\in\cM$. 
Furthermore, it is called nonnegative if $P^i\geq0, \ t\in[0,T]$, a.s.; and called uniformly positive if $P^i\geq c, \ t\in[0,T]$, a.s. with some deterministic constant $c>0$.
\end{definition}

Let
\begin{align*}
&\widetilde A^i=A^i-\lambda E^i-\frac{1}{2}\lambda, \ \widetilde B^i=B^i-\lambda F^i, \ \widetilde Q^i=Q^i+\lambda G^{a,i}(E^i+1)^2,\\
& \widetilde R^i=R^i+\lambda G^{a,i}F^i(F^i)', \ S^i=\lambda G^{a,i}(E^i+1)F^i.
\end{align*}
Then \eqref{PBM} can be written as
\begin{align*}
\begin{cases}
P^{b,i}_t
=G^{b,i}+\int_{t}^{T}\Big[(2\widetilde A^i_s+(C^i_s)'C^i_s)P^{b,i}_s
+2(C^i_s)'\Lambda^{b,i}_s+\widetilde Q^i_s+\sum\limits_{j=1}^\ell q_{ij}P^{b,j}_s\nonumber\\
\qquad\qquad\qquad\qquad-\widehat{\mathcal{N}}(s,P^{b,i}_s,\Lambda^{b,i}_s,i)'\widehat{\mathcal{R}}(s,P^{b,i}_s,i)^{-1}
\widehat{\mathcal{N}}(s,P^{b,i}_s,\Lambda^{b,i}_s,i)\Big]ds
-\int_{t}^{T}(\Lambda^{b,i}_s)'dW_s,\\
\widehat{\mathcal{R}}(s,P^{b,i}_s,i)>0, \ \mbox{for $s\leq T$ and all $i\in\cM$},
\end{cases}
\end{align*}
where
\begin{align*}
\widehat{\mathcal{N}}(s,P,\Lambda,i)&=P((D^i_s)' C^i_s+\widetilde B^i_s)+(D^i_s)'\Lambda
+ S^i_s,\\
\widehat{\mathcal{R}}(s,P,i)&=\widetilde R^i_s+P(D^i_s)'D^i_s.
\end{align*}
And it is associated with the following LQ stochastic control problem on $[0,T]$:
\begin{align*}
\begin{cases}
\mathrm{Minimize} &\ \widetilde J(x, i_0, u)\\
\mbox{subject to} &\ (X, u) \mbox{ is an admissible pair for} \ \eqref{statebm},
\end{cases}
\end{align*}
where the state process is
\begin{align}
\begin{cases}
dX_t=\left[\widetilde A^{\alpha_t}_tX_t+u_t'\widetilde B^{\alpha_t}_t\right]dt
+\left[C^{\alpha_t}_tX_t+D^{\alpha_t}_tu_t\right]'dW_t, \ t\in[0, T], \\
X_0=x, \ \alpha_0=i_0,
\end{cases}
\label{statebm}
\end{align}
and the cost functional is given as the following quadratic form
\begin{align}\label{costbm}
\widetilde J(x, i_0, u(\cdot)):=\mathbb{E}\left\{\int_0^{T}\Big(\widetilde Q_t^{\alpha_t}X_t^2+u_t'\widetilde R_t^{\alpha_t}u_t+2u_t' S^{\alpha_t}_tX_t\Big)dt
+G^{b,\alpha_{T}}X_{T}^2\right\}.
\end{align}

Compared with the stochastic LQ problem studied in our previous work \cite{HSX}, the cross term $2u_t' S^{\alpha_t}_tX_t$ involves in \eqref{costbm}.
The existence and uniqueness of solution to \eqref{PBM} are proved in Theorems 3.5 and 3.6 of \cite{HSX} when $S^i\equiv0$. And the method could be applied to the case
of $S^i=\lambda G^{a,i}(E^i+1)F^i$, however we need to carefully show that the solution is nonnegative or uniformly positive in different cases.

\begin{theorem}[Standard Case]\label{standard}
Under Assumption \ref{assu3}, the system of BSDEs \eqref{PBM} admits a unique nonnegative solution $(P^{b,i},\Lambda^{b,i})_{i\in\cM}$.
\end{theorem}
\pf
For $i\in\cM$, $t\in[0,T]$, $P\in\mathbb{R}^{\ell}$, and $\Lambda\in\mathbb{R}^{n\times \ell}$, set
\begin{align*}
\overline f(t,P, \Lambda, i)=(2\widetilde A^i_t+(C^i_t)'C^i_t+q_{ii})P^i+2(C^i_t)'\Lambda^i+\widetilde Q^i_t+\sum_{j\neq i}q_{ij}P^j.
\end{align*}
From Theorem 3.5 of \cite{HSX}, there exists a unique solution $(\overline P^i, \overline\Lambda^i)_{i\in\cM}$ to the corresponding $\ell$-dimensional linear BSDE with the generator $\overline f$ and terminal value $G^b$, and there exists a constant $c>0$ such that $\overline P^i\leq c$.
Hereafter, we shall use $c$ to represent a generic positive constant independent of $i$, $n$ and $t$, which can be different from line to line.

Denote
\[
H(t,P,\Lambda,i)=-\widehat{\mathcal{N}}(t,P,\Lambda,i)'\widehat{\mathcal{R}}(t,P,i)\widehat{\mathcal{N}}(t,P,\Lambda,i).
\]
For $k\geq1$, $(t, P, \Lambda, i)\in[0, T]\times\mathbb R\times\mathbb{R}^n\times\cM$, define
\[
H^{k}(t, P, \Lambda, i)=\sup_{\widetilde P\in\mathbb R, \; \widetilde\Lambda\in\mathbb R^n}\Big\{H (t, \widetilde P, \widetilde\Lambda, i)-k|P-\widetilde P|-k|\Lambda-\widetilde\Lambda|\Big\}.
\]
Then
it is non-positive and uniformly Lipschitz in $(P, \Lambda)$, and decreases to $H(t, P, \Lambda, i)$ as $k$ goes to infinity.

The following BSDE
\begin{align*}
\begin{cases}
dP^{k,i}_t=-\Big[
\overline f(t,P^k_t, \Lambda^k_t, i)+H^{k}(t,P^{k,i}_t, \Lambda^{k,i}_t, i)\Big]dt+(\Lambda^{k,i}_t)'dW_t, \\
P^{k,i}_T=G^{b,i}, \ \mbox{ for all $i\in\cM$, }
\end{cases}
\end{align*}
is an $\ell$-dimensional BSDE with a Lipschitz generator, so it admits a unique solution, denoted by $\big(P^{k,i}, \Lambda^{k,i}\big)_{i\in\cM}$.
Notice that
\begin{align*}
H^k(t, 0, 0, i)\geq H (t, 0, 0, i),
\end{align*}
and thanks to Lemma \ref{A1} and Assumption \ref{assu3},
\begin{align*}
1-\lambda_t G^{a,i}_t(F^i_t)'[R^i_t+\lambda_t G^{a,i}_tF^i_t(F^i_t)']^{-1}F^i_t=\frac{\mathrm{det}(R^i_t)}{\mathrm{det}(R^i_t+\lambda_t G^{a,i}_tF^i_t(F^i_t)')}>0.
\end{align*}
Hence, we have the following estimate
\begin{align*}
&\quad\;\overline f(t,0,0,i)+H^k(t, 0, 0, i)\\
&\geq\widetilde Q^i_t-\widehat{\mathcal{N}}(t,0,0,i)'\widehat{\mathcal{R}}(t,0,i)^{-1}
\widehat{\mathcal{N}}(t,0,0,i)\\
&=Q^i_t+\lambda_t G^{a,i}_t(E^i_t+1)^2-[\lambda_t G^{a,i}_t(E^i_t+1)F^i_t]'[R^i_t+\lambda_t G^{a,i}_tF^i_t(F^i_t)']^{-1}[\lambda_t G^{a,i}_t(E^i_t+1)F^i_t]\\
&=Q^i_t+\lambda_t G^{a,i}_t(E^i_t+1)^2-(\lambda_t G^{a,i}_t)^2(E^i_t+1)^2(F^i_t)'[R^i_t+\lambda_t G^{a,i}_tF^i_t(F^i_t)']^{-1}F^i_t\\
&=Q^i_t+\lambda_t G^{a,i}_t(E^i_t+1)^2\Big[1-\lambda_t G^{a,i}_t(F^i_t)'[R^i_t+\lambda_t G^{a,i}_tF^i_t(F^i_t)']^{-1}F^i_t\Big]\\
&=Q^i_t+\lambda_t G^{a,i}_t(E^i_t+1)^2\frac{\mathrm{det}(R^i_t)}{\mathrm{det}(R^i_t+\lambda_t G^{a,i}_tF^i_t(F^i_t)')}\\
&\geq0.
\end{align*}

Now we have
$\overline f(t,0,0,i)+H^k(t, 0, 0, i)\geq0, \ G\geq 0$, and
\[
\overline f(t, P, \Lambda, i)+H^k(t,P^i, \Lambda^i, i) \leq \overline f(t, P, \Lambda, i),
\]
so by the comparison theorem for multidimensional BSDEs (see Lemma 3.4 of \cite{HSX} or Lemma 2.2 of \cite{HLT}), we have
\begin{align*}
\label{Pbound}
0\leq P^{k,i}_t\leq\overline P^i_t\leq c,
\end{align*}
and also $P^{k,i}_t$ is decreasing in $k$, for each $i\in\cM$.

Let $P^i_t=\lim\limits_{k\rightarrow\infty}P^{k,i}_t$, $i\in\cM$.
It is important to note that we can regard $\big(P^{k,i}, \Lambda^{k,i}\big)$ as the solution of a scalar-valued quadratic BSDE for each fixed $i\in\cM$. Thus by Proposition 2.4 in \cite{Ko}, there exists a process $\Lambda\in L^{2}_{\mathbb F}(0, T;\mathbb{R}^{n\times \ell})$ such that $(P, \Lambda)$ is a solution to BSDE \eqref{PBM}. We have now established the existence of the solution.

The uniqueness part is similar to Theorem 3.5 of \cite{HSX}, so we omit it.
\eof

\begin{theorem}[Singular Case I]\label{singular}
Under Assumption \ref{assu2}, the system of BSDEs \eqref{PBM} admits a unique uniformly positive solution $(P^{b,i},\Lambda^{b,i})_{i\in\cM}$.
\end{theorem}
\pf
Consider the following $\ell$-dimensional decoupled BSDEs:
\begin{equation}
\begin{cases}
dP^{b,i}=-\Big[(2\widetilde A^i+(C^i)'C^i+q_{ii})P^{b,i}+2(C^i)'\Lambda^{b,i}+\widetilde Q^i+H(P^{b,i}, \Lambda^{b,i}, i)\Big]dt+(\Lambda^{b,i})'dW,\\
P^{b,i}_T=G^{b,i},\\
P^{b,i}>0, \ \mbox{ for all $i\in\cM$}.
\end{cases}
\label{P:lower}
\end{equation}
If we can show that \eqref{P:lower} admits a uniformly positive solution $(\underline P^i, \underline\Lambda^i)_{i\in\cM}$, then a solution to \eqref{PBM} can be constructed following the procedure of Theorem 3.6 in \cite{HSX}. So it remains to study the solvability of \eqref{P:lower}.

As the system of BSDEs \eqref{P:lower} is decoupled, we would only consider the solvability of each fixed one of the $\ell$ equations in \eqref{P:lower} and may drop the superscript $i$ in the remaining proof for notation simplicity.

Let us first consider the following BSDE:
\begin{align}\label{Psin2}
\begin{cases}
dP_t=-\Big[- [P(D' C+\widetilde B)+D'\Lambda+ S]'[\lambda G^aFF'+PD'D]^{-1}[P(D' C+\widetilde B)+D'\Lambda+ S]\\
\qquad\qquad\qquad+(2\widetilde A+C'C+q_{ii})P
+2C'\Lambda+\widetilde Q\Big]dt+\Lambda'dW,\\
P_T=G^b,\\
P>0.
\end{cases}
\end{align}
Under Assumption \ref{assu2}, $D'D$ is invertible, so the generator $f$ of \eqref{Psin2} can be rewritten as
\begin{align*}
f(P,\Lambda)&=(2\widetilde A+C'C+q_{ii})P
+2C'\Lambda+\widetilde Q\\
&\qquad- [P(D' C+\widetilde B)+D'\Lambda+ S]'[\lambda G^aFF'+PD'D]^{-1}[P(D' C+\widetilde B)+D'\Lambda+ S]\\
&=(2\widetilde A+C'C+q_{ii})P
+2C'\Lambda+\widetilde Q-P(D'C+\widetilde B)'(D'D)^{-1}(D'C+\widetilde B)\\
&\qquad-2(D'C+\widetilde B)'(D'D)^{-1}(S+D'\Lambda)+\lambda G^a|(D'C+\widetilde B)'(D'D)^{-1}F|^2\\
&\qquad-[S-\lambda G^aFF'(D'D)^{-1}(D'C+\widetilde B)+D'\Lambda]'(PD'D+\lambda G^aFF')^{-1}\\
&\qquad\qquad \times[S-\lambda G^aFF'(D'D)^{-1}(D'C+\widetilde B)+D'\Lambda],
\end{align*}
by some basis operations of positive matrices.

Under Assumption \ref{assu2}, there is a deterministic constant $c_2>0$ such that
\begin{align*}
|2 A-2\lambda E-\lambda+C'C+q_{ii}- (D'C+ B-\lambda F)'(D'D)^{-1}(D'C+B-\lambda F)|\leq c_2.
\end{align*}
Let $\varepsilon=\delta e^{-c_2T}$, where $\delta>0$ is the constant in Assumption \ref{assu2}.
Then from Theorem 2.3 of \cite{Ko}, there is a bounded maximal solution (see page 565 of \cite{Ko} for its definition) $(P^{\varepsilon},\Lambda^{\varepsilon})$ to the BSDE with terminal value $G^b$ and generator $f^{\varepsilon}(P,\Lambda)$, where
\begin{align*}
f^{\varepsilon}(P,\Lambda)&=(2\widetilde A+C'C+q_{ii})P
+2C'\Lambda+\widetilde Q- P(D'C+\widetilde B)'(D'D)^{-1}(D'C+\widetilde B)\\
&\qquad-2(D'C+\widetilde B)'(D'D)^{-1}(S+D'\Lambda)+\lambda G^a|(D'C+\widetilde B)'(D'D)^{-1}F|^2\\
&\qquad-[S-\lambda G^aFF'(D'D)^{-1}(D'C+\widetilde B)+D'\Lambda]'((P\vee\varepsilon)D'D+\lambda G^aFF')^{-1}\\
&\qquad\qquad \times[S-\lambda G^aFF'(D'D)^{-1}(D'C+\widetilde B)+D'\Lambda].
\end{align*}
Notice that $Q\geq0$, and
\begin{align*}
&\quad\;\lambda G^a(E+1)^2-2\lambda G^a(E+1)(D'C+B-\lambda F)'(D'D)^{-1}F+\lambda G^a|(D'C+B-\lambda F)'(D'D)^{-1}F|^2\\
&\qquad-\lambda^2 (G^a)^2[(E+1)F-FF'(D'D)^{-1}(D'C+B-\lambda F)]'((P\vee\varepsilon)D'D+\lambda G^aFF')^{-1}\\
&\qquad\qquad \times[(E+1)F-FF'(D'D)^{-1}(D'C+B-\lambda F)]\\
&=\lambda G^a|E+1-(D'C+B-\lambda F)'(D'D)^{-1}F|^2\\
&\qquad-\lambda^2 (G^a)^2[E+1-F'(D'D)^{-1}(D'C+B-\lambda F)]'F'((P\vee\varepsilon)D'D+\lambda G^aFF')^{-1}F\\
&\qquad\qquad \times[E+1-F'(D'D)^{-1}(D'C+B-\lambda F)]\\
&=\lambda G^a|E+1-(D'C+B-\lambda F)'(D'D)^{-1}F|^2(1-\lambda G^aF'((P\vee\varepsilon)D'D+\lambda G^aFF')^{-1}F)\\
&=\lambda G^a|E+1-(D'C+B-\lambda F)'(D'D)^{-1}F|^2\frac{\mathrm{det}((P\vee\varepsilon)D'D)}{\mathrm{det}((P\vee\varepsilon)D'D+\lambda G^{a}FF')}\\
&\geq0,
\end{align*}
where the last equality is due to Lemma \ref{A1}, so we have
\begin{align*}
f^{\varepsilon}(P,\Lambda)&=(2 A-2\lambda E-\lambda+C'C+q_{ii})P
+2C'\Lambda+ Q+\lambda G^a(E+1)^2\\
&\qquad- P(D'C+ B-\lambda F)'(D'D)^{-1}(D'C+B-\lambda F)-2\lambda G^a(E+1)(D'C+B-\lambda F)'(D'D)^{-1}F\\
&\qquad+\lambda G^a|(D'C+B-\lambda F)'(D'D)^{-1}F|^2\\
&\qquad-\lambda^2 (G^a)^2[E+1-F'(D'D)^{-1}(D'C+B-\lambda F)]'F'((P\vee\varepsilon)D'D+\lambda G^aFF')^{-1}F\\
&\qquad\qquad \times[E+1-F'(D'D)^{-1}(D'C+B-\lambda F)]\\
&\qquad-2(D'C+ B-\lambda F)'(D'D)^{-1}D'\Lambda-\Lambda'D((P\vee\varepsilon)D'D+\lambda G^aFF')^{-1}D'\Lambda\\
&\qquad-2\lambda G^a[E+1-F'(D'D)^{-1}(D'C+B-\lambda F)]F'((P\vee\varepsilon)D'D+\lambda G^aFF')^{-1}D'\Lambda\\
&\geq(2 A-2\lambda E-\lambda+C'C+q_{ii})P- (D'C+ B-\lambda F)'(D'D)^{-1}(D'C+B-\lambda F)P
+2C'\Lambda\\
&\qquad-2(D'C+ B-\lambda F)'(D'D)^{-1}D'\Lambda-\Lambda'D((P\vee\varepsilon)D'D+\lambda G^aFF')^{-1}D'\Lambda\\
&\qquad-2\lambda G^a[E+1-F'(D'D)^{-1}(D'C+B-\lambda F)]F'((P\vee\varepsilon)D'D+\lambda G^aFF')^{-1}D'\Lambda.
\end{align*}

Obviously, the following BSDE
\begin{align*}
\begin{cases}
dP=-\Big[-c_2P+2C'\Lambda-2(D'C+ B-\lambda F)'(D'D)^{-1}D'\Lambda-\Lambda'D((P\vee\varepsilon)D'D+\lambda G^aFF')^{-1}D'\Lambda\\
\qquad\qquad\quad-2\lambda G^a[E+1-F'(D'D)^{-1}(D'C+B-\lambda F)]F'((P\vee\varepsilon)D'D+\lambda G^aFF')^{-1}D'\Lambda\Big]dt+\Lambda'dW,\\
P_T=\delta
\end{cases}
\end{align*}
admits a solution $(\delta e^{-c_2(T-t)},0)$. Hence the maximal solution argument in Theorem 2.3 of \cite{Ko} gives
\begin{align}\label{aa}
P^{\varepsilon}_t\geq \delta e^{-c_2(T-t)}\geq \delta e^{-c_2T}=\varepsilon.
\end{align}
This implies that $(P^{\varepsilon},\Lambda^{\varepsilon})$ is actually a solution to \eqref{Psin2}.

On the other hand,
\begin{align*}
- [P(D' C+\widetilde B)+D'\Lambda+ S]'[\lambda G^aFF'+PD'D]^{-1}[P(D' C+\widetilde B)+D'\Lambda+ S]\leq 0,
\end{align*}
and the following linear BSDE
\begin{align}\label{BSDElinear}
\begin{cases}
dP_t=-\Big[(2\widetilde A+C'C+q_{ii})P+2C'\Lambda+\widetilde Q\Big]dt+\Lambda'dW,\\
P_T=G^b,
\end{cases}
\end{align}
admits a unique solution solution $(\overline P, \overline\Lambda)$ such that $\overline P\leq c_3$ for some constant $c_3>0$. Then applying the comparison theorem to BSDEs {\eqref{Psin2} and \eqref{BSDElinear}}, we get
\begin{align*}
P^{\varepsilon}_t\leq \overline P_{t}\leq c_3.
\end{align*}
This implies that $(P^{\varepsilon},\Lambda^{\varepsilon})$ is actually a solution to the following BSDE:
\begin{align}
\label{Psin3}
\begin{cases}
dP_t=-\Big[- [P(D' C+\widetilde B)+D'\Lambda+ S]'[\lambda G^aFF'+PD'D]^{-1}[P(D' C+\widetilde B)+D'\Lambda+ S]g^{\varepsilon,c_3}(P)\\
\qquad\qquad\quad+(2\widetilde A+C'C+q_{ii})P
+2C'\Lambda+\widetilde Q\Big]dt+\Lambda'dW,\\
P_T=G^b,
\end{cases}
\end{align}
where $g^{\varepsilon,c_3}:\mathbb{R}^+\rightarrow[0,1]$ is a smooth truncation function
satisfying $g^{\varepsilon,c_3}(x)=1$ for $x\in[\varepsilon,c_3]$, and $g^{\varepsilon,c_3}(x)=0$ for $x\in[0,\varepsilon/2]\cup[2c_3,\infty)$.

Finally, under Assumption \ref{assu2}, we have
\begin{align*}
&\quad[P(D' C+\widetilde B)+D'\Lambda+ S]'[R+\lambda G^aFF'+PD'D]^{-1}[P(D' C+\widetilde B)+D'\Lambda+ S]g^{\varepsilon,c_3}(P)\\
&\leq \frac{2\delta}{\varepsilon}|P(D' C+\widetilde B)+D'\Lambda+ S|^2g^{\varepsilon,c_3}(P)\\
&\leq c(1+|\Lambda|^2).
\end{align*}
By Theorem 2.3 in \cite{Ko}, there exists a bounded, maximal solution $(P_2,\Lambda_2)$ to the following quadratic BSDE:
\begin{align}
\label{Psin4}
\begin{cases}
dP_t=-\Big[- [P(D' C+\widetilde B)+D'\Lambda+ S]'[R+\lambda G^aFF'+PD'D]^{-1}[P(D' C+\widetilde B)+D'\Lambda+ S]g^{\varepsilon,c_3}(P)\\
\qquad\qquad\quad+(2\widetilde A+C'C+q_{ii})P
+2C'\Lambda+\widetilde Q\Big]dt+\Lambda'dW,\\
P_T=G^b.
\end{cases}
\end{align}
Notice that under Assumption \ref{assu2}, $R\geq0$. Applying the maximal solution argument in Theorem 2.3 of \cite{Ko} to BSDEs \eqref{Psin3} and \eqref{Psin4}, we get
\begin{align}\label{bb}
P^{\varepsilon}\leq P_2.
\end{align}
On the other hand,
\begin{align*}
- [P(D' C+\widetilde B)+D'\Lambda+ S]'[R+\lambda G^aFF'+PD'D]^{-1}[P(D' C+\widetilde B)+D'\Lambda+ S]g^{\varepsilon,c_3}(P)\leq 0,
\end{align*}
thus applying the comparison theorem to BSDEs \eqref{BSDElinear} and \eqref{Psin4}, we have
\begin{align}\label{cc}
P_2\leq \overline P_{t}\leq c_3.
\end{align}
Combining \eqref{aa}, \eqref{bb}, \eqref{cc}, we obtain $g^{\varepsilon,c_3}(P_{2})=1$. Hence $(P_2,\Lambda_2)$ is actually a solution to \eqref{P:lower}.

The uniqueness part is similar to Theorem 3.6 of \cite{HSX}, so we omit it.
\eof

\begin{theorem}[Singular Case II]\label{singuspec}
Under Assumption \ref{assu4}, the system of BSDEs \eqref{PBM} admits a unique nonnegative solution
$(P^{b,i},\Lambda^{b,i})_{i\in\cM}$.
\end{theorem}
\pf
Notice that under Assumption \ref{assu4},
\begin{align*}
1-\frac{\lambda G^{a,i}(F^i)^2}{R^i+\lambda G^{a,i}(F^i)^2}=\frac{R^i}{R^i+\lambda G^{a,i}(F^i)^2}\geq0,
\end{align*}
then by Theorem \ref{standard}, BSDE \eqref{PBM} admits a unique nonnegative solution.
\eof

The following remark will be used in Section 5.
\begin{remark}\label{sinuni}
If Assumption \ref{assu4} holds with $G^{i,b}\geq \delta>0$, then the solution
$(P^{b,i},\Lambda^{b,i})_{i\in\cM}$ of the system of BSDEs \eqref{PBM} is uniformly positive.

The argument is similar to the proof of Theorem \ref{singular}, so we will only give the sketch showing that solution to \eqref{PBM} has a positive lower bound given by the following BSDE (the superscript $i$ is suppressed):
\begin{align}\label{Psinspec}
\begin{cases}
dP_t=-\Big[- \frac{1}{\lambda G^aF^2}[P(D'C+ B-\lambda F)+D'\Lambda+ \lambda G^{a}(E+1)F]^2\\
\qquad\qquad+(2A-2\lambda E-\lambda+|C|^2+q_{ii})P
+2C'\Lambda+ Q+\lambda G^{a}(E+1)^2\Big]dt+\Lambda' dW,\\
P_T=G^b,\\
P>0.
\end{cases}
\end{align}
By the proof of Theorem \ref{standard}, the truncation method and Theorem 2.3 in \cite{Ko}, there exists a bounded, maximal solution $(P,\Lambda)$ to \eqref{Psinspec}, such that $P\leq c$ for some deterministic constant $c>0$.

Under Assumption \ref{assu4}, there exists a constant $c_4>0$ such that
\[
|2A-2\lambda E-\lambda+C^2+q_{ii}-\frac{2}{F}(D'C+ B-\lambda F)(E+1)|\leq c_4,
\]
and
\[
\frac{1}{\lambda G^aF^2}(D'C+ B-\lambda F)^2\leq c_4.
\]
Hence, for $P>0$ and $\Lambda\in\mathbb{R}^n$, we have the following estimate for the generator of \eqref{Psinspec}:
\begin{align*}
g(P,\Lambda):&=- \frac{1}{\lambda G^aF^2}[P(D'C+ B-\lambda F)+D'\Lambda+ \lambda G^{a}(E+1)F]^2\\
&\qquad\qquad+(2A-2\lambda E-\lambda+|C|^2+q_{ii})P
+2C'\Lambda+ Q+\lambda G^{a}(E+1)^2\\
&\geq (2A-2\lambda E-\lambda+|C|^2+q_{ii}-\frac{2}{F}(D'C+ B-\lambda F)(E+1))P-\frac{1}{\lambda G^aF^2}(D'C+ B-\lambda F)^2P^2\\
&\qquad\qquad+2C'\Lambda-\frac{2}{F}(E+1)D'\Lambda-\frac{2}{\lambda G^aF^2}P(D'C+ B-\lambda F)D'\Lambda-\frac{1}{\lambda G^aF^2}|D'\Lambda|^2\\
&\geq-c_4P-c_4P^2+2C'\Lambda-\frac{2}{F}(E+1)D'\Lambda-\frac{2}{\lambda G^aF^2}P(D'C+ B-\lambda F)D'\Lambda-\frac{1}{\lambda G^aF^2}|D'\Lambda|^2.
\end{align*}

Obviously, the following BSDE
\begin{align*}
\begin{cases}
dP=-\Big[-c_4P-c_4P^2+2C'\Lambda-\frac{2}{F}(E+1)D'\Lambda-\frac{2}{\lambda G^aF^2}(DC+ B-\lambda F)PD'\Lambda-\frac{1}{\lambda G^aF^2}|D'\Lambda|^2\Big]dt+\Lambda'dW,\\
P_T=\delta
\end{cases}
\end{align*}
admits a solution
$(\frac{1}{(1+\frac{1}{\delta})e^{c_4(T-t)}-1},0)$. Hence the maximal solution argument in Theorem 2.3 of \cite{Ko} gives
\[
P\geq \frac{1}{(1+\frac{1}{\delta})e^{c_4(T-t)}-1}\geq \frac{1}{(1+\frac{1}{\delta})e^{c_4T}-1}>0.
\]
\end{remark}

Combining Lemma \ref{link}, Theorems \ref{standard}, \ref{singular} and \ref{singuspec}, we have the following solvability of the system of BSDEs \eqref{P}.
\begin{theorem}\label{existence}
Under Assumption \ref{assu3} (resp. \ref{assu2}, \ref{assu4}), let $(P^{b,i},\Lambda^{b,i})_{i\in\cM}$ be the unique nonnegative (resp. uniformly positive, nonnegative) solution of \eqref{PBM} and define
\begin{align*}
&P^{i}_t:=P^{b,i}_t\mathbf{1}_{t<\tau}+G^{a,i}_{\tau}\mathbf{1}_{t\geq\tau},\\
&\Lambda^{i}_t:=\Lambda^{b,i}_t\mathbf{1}_{t\leq\tau},\\
&U^{i}_t:=(G^{a,i}_t-P^{b,i}_t)\mathbf{1}_{t\leq\tau}, \ t\in[0,T], \ i\in\cM.
\end{align*}
Then $(P^{i},\Lambda^{i},U^{i})_{i\in\cM}$ is a nonnegative (resp. uniformly positive, positive) solution of \eqref{P}.
\end{theorem}
\pf
It is obvious $(P^i, \Lambda^i, U^i)\in S^\infty_{\mathbb{G}}(0, T; \mathbb {R})\times L^{2}_{\mathbb{G}}(0, T;\mathbb{R}^{n})\times L^{\infty}_{\mathbb{G}}(0,T;\mathbb{R}))$, for all $i\in\cM$.
By Lemma \ref{link}, $(P^{i},\Lambda^{i},U^{i})_{i\in\cM}$ satisfies the first equality of \eqref{P}.
It remains to show $\mathcal{R}(t,P^i_{t-},U^i_t, i)>0$, for $t\leq T\wedge\tau$, $i\in\cM$.

By the definition of $\mathcal{R}(t,P^i_{t-},U^i_t, i)$, and $P^{i}_{t-}+U^{i}_t=G^{a,i}_t \mathbf{1}_{t\leq\tau}+G^{a,i}_{\tau}\mathbf{1}_{t>\tau}=G^{a,i}_{\tau\wedge t}$, we have
\begin{align*}
\mathcal{R}(t,P^i_{t-},U^i_t,i)&=R^i_t+P^i_{t-}(D^i_t)'D^i_t+\lambda^{\mathbb{G}}_t(P^i_{t-}+U^i_t)(F^i_t)(F^i_t)'\\
&=R^i_t+P^i_{t-}(D^i_t)'D^i_t+\lambda^{\mathbb{G}}_tG^{a,i}_{\tau\wedge t}(F^i_t)(F^i_t)'.
\end{align*}
In Standard Case, we have $R^i_t\geq\delta I_{m}$, $G^i_t\geq 0$ and $P^i_t\geq0$, thus
$\mathcal{R}(t,P^i_{t-},U^i_t, i)>0, \ i\in\cM$.
In Singular Case I, we have $R^i_t\geq0$, $G^i\geq \delta$ , $(D^i_t)'D^i_t\geq\delta I_{m}$ and $P^i_t\geq c>0$, thus $\mathcal{R}(t,P^i_{t-},U^i_t, i)>0, \ i\in\cM$.
In Singular Case II, we have $R^i_t\geq0$, $G^{i,a}_t\geq \delta$, $\lambda^{\mathbb{G}}_t(F^i_t)^2=\lambda_t(F^i_t)^2\geq\delta$ and $P^i_t\geq 0$ for $t\leq T\wedge\tau$, thus $\mathcal{R}(t,P^i_{t-},U^i_t, i)>0$, for $t\leq T\wedge\tau$, $i\in\cM$.
\eof
\begin{remark}
A matrix-valued SRE with Poisson jumps was solved in Zhang, Dong and Meng \cite{ZDM} with uniformly positive control weight (corresponding to our Standard Case), while a scalar-valued SRE with a jump and $m=n=1$, $C=0,\ E=Q=0, \ R=0$ was studied in Kharroubi, Lim and Ngoupeyou \cite{KLN} (corresponding to our Singular Case I). Up to our knowledge, no existing literature has concerned the solvability of SRE or stochastic LQ problem with jumps corresponding to our Singular Case II.

Recall the definition of $\mathcal{R}(t,P,U,i)$ in \eqref{R}, $$\mathcal{R}(t,P,U,i)=R^i+P(D^i)'D^i+\lambda^{\mathbb{G}}(P+U)F^i(F^i)'.$$
There are three components $R^i$, $P(D^i)'D^i$ and $\lambda^{\mathbb{G}}(P+U)F^i(F^i)'$ on the right hand. By Theorem \ref{existence}, Assumptions \ref{assu3}, \ref{assu2} and \ref{assu4} will result in the uniformly positivity of one of these three components, while nonnegativity of the other two terms, hence the second constraint $\mathcal{R}(t,P,U, i)>0$ in \eqref{P} is fulfilled.
\end{remark}

\section{Solution to the stochastic LQ problem \eqref{LQ}}
We slightly strengthen Assumption \ref{assu4} to the following:
\begin{assumption}[Singular Case II']
\label{assu5}
$m=1$, $Q^i\geq0$, $R^i\geq0$, and
there exists a constant $\delta>0$ such that $G^{i} \geq \delta$ and $\lambda(F^i)^2\geq\delta$, for all $i\in\cM$.
\end{assumption}

\begin{theorem}\label{LQcontrol}
Under Assumption \ref{assu3} (resp. \ref{assu2} and \ref{assu5}), let $(P^i, \Lambda^i, U^i)_{i\in\cM}$ be nonnegative (resp. uniformly positive, positive) solution of \eqref{P}. Then the problem \eqref{LQ} has an optimal control, as a feedback function of the time $t$, the state $X$ and the market regime $i$,
\begin{align}
\label{opticon}
u^*(t, X_{t-}, i)&=-\mathcal{R}(t,P^i_{t-},U^i_t,i)^{-1}\mathcal{N}(t,P^i_{t-},\Lambda^i_t,U^i_t,i)X_{t-}.
\end{align}
Moreover, the corresponding optimal value is
\begin{align}\label{optival}
V(x,i_0)=\min_{u\in\mathcal{U}}J(x, i_0, u)=P_0^{i_0}x^2.
\end{align}
\end{theorem}
\pf
In light of the length of many equations, $``t,X_{t-},\alpha_t,P^{\alpha_t}_{t-},\Lambda^{\alpha_t},U^{\alpha_t}"$ will be suppressed where no confusion occurs in the sequel.

Substituting \eqref{opticon} into the state process \eqref{state} (with $``i"$ replaced by $``\alpha_t"$), we have
\begin{align}
\label{statefeed}
\begin{cases}
dX_s^*=\left[A^{\alpha_s}_s-(B^{\alpha_s}_s)'\mathcal{R}^{-1}\mathcal{N}\right]X^*_{s-} ds+\left[C^{\alpha_s}_s-D^{\alpha_s}_s\mathcal{R}^{-1}\mathcal{N}\right]'X^*_{s-}dW_s\\
\qquad\qquad\qquad+\left[E^{\alpha_s}_s-(F^{\alpha_s}_s)'\mathcal{R}^{-1}\mathcal{N}\right]X^*_{s-}dM_s, \ s\in[0,T\wedge\tau], \\
X^*_0=x, \ \alpha_0=i_0.
\end{cases}
\end{align}
Since the coefficients of SDE \eqref{statefeed} are square integrable w.r.t. $t$, it admits a unique strong solution $X^*$.
Actually
\begin{align*}
\begin{cases}
X_t^*=x\Phi_t, \ t\in[0,T\wedge\tau),\\
X_{T\wedge\tau}^*=x\Phi_T\mathbf{1}_{\tau> T} +x\Phi_{\tau}(1+E_\tau-F_\tau'\mathcal{R}^{-1}\mathcal{N})\mathbf{1}_{\tau\leq T},
\end{cases}
\end{align*}
is the solution to \eqref{statefeed},
where
\begin{align*}
\Phi_t=\exp\Big(\int_0^t\Big(A-\lambda^{\mathbb{G}}E+(\lambda^{\mathbb{G}}F-B)'\mathcal{R}^{-1}\mathcal{N}
-\frac{1}{2}|C-D\mathcal{R}^{-1}\mathcal{N}|^2\Big)ds+\int_0^t(C-D\mathcal{R}^{-1}\mathcal{N})'dW_s\Big).
\end{align*}

Let $u^*_{t}=u^*(t,X^{*}_{t-},\alpha_t)$.
Define a sequence of stopping times $\{\iota_j\}_{j\geq1}$ as follows:
\begin{align*}
\iota_j:=\inf\Big\{t>0: |X^*_{t-}|+|u_{t }^*|> j\Big\}\wedge j,
\end{align*}
with the convention that $\inf\emptyset=\infty$. It is obvious that $\iota_j\uparrow\infty$, a.s. as $j\uparrow\infty$. Applying It\^{o}'s formula to $P^{\alpha_t}_t (X^*_t)^2$, we have
for any stopping time $\iota\leq T\wedge\tau$,

\begin{align}
\label{bound}
P_0^{i_0}x^2=\E\Big[P^{\alpha_{\iota\wedge\iota_j}}_{\iota\wedge\iota_j}(X_{\iota\wedge\iota_j}^*)^2
+\int_0^{\iota\wedge\iota_j}\Big(Q_s(X_{s-}^*)^2+(u_s^*)'R_su_s^*\Big)ds\Big].
\end{align}

In Standard Case, $Q^i\geq0, \ R^i\geq\delta I_{m}$, it follows
\begin{align*}
\delta\E\int_0^{T\wedge\tau\wedge\iota_j}|u^*_{t}|^2dt
\leq \E\int_0^{T\wedge\tau\wedge\iota_j}(u_s^*)'R_su_s^*dt
\leq P_0^{i_0}x^2.
\end{align*}
Sending $j\uparrow\infty$, by the monotone convergence theorem, we have $u^*_{t}\in L^2_{\mathbb{H}}(0,T\wedge\tau;\mathbb{R}^m)$.

In Singular Case I (resp. Singular Case II'), there exists some constant $c>0$ such that $P_t^i\geq c$, for all $i\in\cM$ by Theorem \ref{singular} (resp. Theorem \ref{singuspec} and Remark \ref{sinuni}). Then from \eqref{bound}, we have
\begin{align*}
c\E\Big[|X_{\iota\wedge\iota_j}^*|^2\Big]\leq P_0^{i_0}x^2.
\end{align*}
By Fatou's lemma, we have
\begin{align*}
\E[|X_{\iota}^*|^2]\leq P_0^{i_0}x^2,
\end{align*}
for any stopping time $\iota\leq T\wedge\tau$. This further implies
\begin{align}\label{intx}
\E\int_0^{T\wedge\tau}|X^*_{t-}|^2dt\leq\int_0^{T\wedge\tau}\E[|X^*_{t-}|^2]dt\leq cT.
\end{align}
Applying It\^{o}'s formula to $|X^*_t|^2$ and using the above two inequalities, we have
\begin{align*}
&\qquad x^2+\E\int_0^{T\wedge\tau\wedge\iota_j}\Big[(u^*)'(D'D+\lambda^{\mathbb{G}}FF')u^*\Big]ds\\
&=\E|X_{T\wedge\tau\wedge\iota_j}^*|^2-\E\int_0^{T\wedge\tau\wedge\iota_j}
\Big[(2A+C'C+\lambda^{\mathbb{G}}E^2)(X_{s-}^{*})^2+ 2(B'+C'D+\lambda^{\mathbb{G}}EF')u^*X^*_{s-} \Big] ds\\
&\leq c+\E\int_0^{T\wedge\tau\wedge\iota_j} 2c|u^*_{s}X^*_{s-} | ds
\end{align*}
By Assumption \ref{assu2} or \ref{assu5}, the elementary inequality $2cab\leq\frac{2c^{2}}{\delta}a^2+\frac{\delta}{2}b^2$ and \eqref{intx}, we have
\begin{align*}
\delta\E\int_0^{T\wedge\tau\wedge\iota_j}|u^*|^2ds &\leq c+ \frac{2c^{2}}{\delta}\E\int_0^{T\wedge\tau\wedge\iota_j} |X_{s-}^*|^2ds +\frac{\delta}{2}\E\int_0^{T\wedge\tau\wedge\iota_j}|u^*|^2ds\\
&\leq c+\frac{\delta}{2}\E\int_0^{T\wedge\tau\wedge\iota_j}|u^*|^2ds.
\end{align*}
After rearrangement and sending $j\uparrow\infty$, it follows from the monotone convergent theorem that $u^*\in L^2_{\mathbb{H}}(0,T\wedge\tau;\mathbb{R}^m)$.

The remaining proof is to verify the optimality of $u^*$ through completion of square technique, it is very similar to Theorem 5.1 in \cite{ZDM} or Theorem 4.2 in \cite{HSX}, so we omit it.
\eof

\section{Mean-variance hedging problem}
We consider a financial market consisting of $m+1$ primitive assets: one risk-free asset with zero interest rate and $m$ risky assets (the stocks) whose price processes follow
\begin{align*}
dS_t=\mathrm{diag}(S_t)\Big(\mu_t^{\alpha_t}dt+\sigma_t^{\alpha_t}dW_t+F_t^{\alpha_t} dM_t\Big), \ t\in[0, T\wedge\tau].
\end{align*}
Assume that the appreciation process $\mu^i\in L^{\infty}_{\mathbb{F}}(0, T;\mathbb{R}^m)$, the volatility $\sigma^i\in L^{\infty}_{\mathbb{F}}(0, T;\mathbb{R}^{m\times n})$ and $F^i\in L_{\mathbb{F}}^\infty(0,T;\mathbb{R}^m)$, $F^i_j\geq -1$ for all $i\in\cM$, $j=1,...,m$. Also there exists a constant $\delta>0$ such that
\begin{itemize}
\item[(i)] $\sigma^i(\sigma^i)'\geq \delta I_{m}$; or
\item[(ii)] $m=1, \ \mbox{and} \ \lambda(F^i)^2\geq \delta$, \ for all \ $i\in\cM$.
\end{itemize}

For any $x\in\R$ and $\pi\in L^{2}_{\mathbb{H}}(0, T\wedge\tau;\mathbb{R}^m)$, the wealth process $X$ for a self-financing investor with initial capital $x$ and with quantity $\pi$ invested in the risky assets is described by
\begin{align}
\label{wealth}
\begin{cases}
dX_t=\pi_t'\mu^{\alpha_t}_tdt+\pi_t'\sigma^{\alpha_t}_tdW_t+\pi_t'F_t^{\alpha_t}dM_t, \ t\in[0, T\wedge\tau],\\
X_0=x, \ \alpha_0=i_0.
\end{cases}
\end{align}

Assume that
for each $i\in\cM$,
$H^i$ is a bounded $\mathcal{G}_{T\wedge\tau}$-measurable random variable of the form
\begin{align}
H^i=H^{b,i}\mathbf{1}_{T<\tau}+H^{a,i}_{\tau}\mathbf{1}_{T\geq\tau},
\end{align}
where
\begin{align}
H^{b,i}\in L^\infty_{\mathcal{F}_T}(\Omega), \ H^{a,i}\in L^\infty_{\mathbb{F}}(0,T;\mathbb{R}).
\end{align}

Consider the following mean-variance hedging problem
\begin{align}
\label{mvhedge}
\min_{\pi\in L^{2}_{\mathbb{H}}(0, T\wedge\tau;\mathbb{R}^m)}\E(X^{\pi}_{T\wedge\tau}-H^{\alpha_{T\wedge\tau}})^2,
\end{align}
where $X^{\pi}$ is the solution to the wealth equation \eqref{wealth}.

Now we apply the general results obtained in the previous section to the mean-variance hedging problem \eqref{mvhedge}, where
\begin{align*}
A=0, \ B=\mu, \ C=0, \ D=\sigma', \ u=\pi, \ E=0, \ Q=0, \ R=0, \ G=1.
\end{align*}

\begin{remark}
If $m=n=1$, and there is no regime switching, the problem \eqref{mvhedge} was solved in \cite{KLN} under the condition (i), namely $\sigma^2\geq\delta>0$.
\end{remark}

The associated system of BSDEs \eqref{P} for the problem \eqref{mvhedge} reduces to the following:
\begin{align}\label{Pmv}
\begin{cases}
P^i_t=1+\int_{t\wedge\tau}^{T\wedge\tau}\Big[-\mathcal{N}(s,P^i_{s-},\Lambda^i_s,U^i_s,i)'
\mathcal{R}(s,P^i_{s-},U^i_s,i)^{-1}
\mathcal{N}(s,P^i_{s-},\Lambda^i_s,U^i_s,i)+\sum_{j=1}^\ell q_{ij}P^j_{s-}\Big]ds\\
\qquad\qquad\qquad\qquad-\int_{t\wedge\tau}^{T\wedge\tau}(\Lambda^i_s)'dW_s-\int_{t\wedge\tau}^{T\wedge\tau}U^i_sdM_s,\\
\mathcal{R}(s,P^i_{s-},U^i_s,i)>0, \ \mbox{for $s\leq T\wedge\tau$ and all $i\in\cM$},
\end{cases}
\end{align}
where
\begin{align*}
\mathcal{N}(s,P,\Lambda,U,i)&=P\mu^i_s+\sigma^i_s\Lambda_s+\lambda^{\mathbb{G}}_sF^i_s U,\\
\mathcal{R}(s,P,U,i)&=P\sigma^i_s(\sigma^i_s)'+\lambda^{\mathbb{G}}_s(P+U)F^i_s(F^i_s)'.
\end{align*}

From Theorems \ref{singular}, \ref{singuspec} and Remark \ref{sinuni}, we know that \eqref{Pmv} admits a uniformly positive solution $(P^i,\Lambda^i,U^i)_{i\in\cM}$, such that
$(P^i, \Lambda^i, U^i)\in S^\infty_{\mathbb{G}}(0, T; \mathbb {R})\times L^{2}_{\mathbb{G}}(0, T;\mathbb{R}^{n})\times L^{\infty}_{\mathbb{G}}(0,T;\mathbb{R}))$ for all $i\in\cM$.

\begin{remark}
For the solution $(P^i,\Lambda^i,U^i)_{i\in\cM}$ of \eqref{Pmv} constructed in Theorem \ref{existence}, we know $P^i_{t-}+U^i_t\equiv 1$, as the terminal value of \eqref{Pmv} is identically equal to $1$. But we will keep $P^i_{t-}+U^i_t$ in the sequel as it is the case for general terminal value of \eqref{Pmv}.
\end{remark}

To construct a solution of the problem \eqref{mvhedge}, we still need to consider the following system of linear BSDEs with jumps:

\begin{align}\label{h}
h^i_t&=H^i-\int_{t\wedge\tau}^{T\wedge\tau}\frac{1}{P^i_{s-}}\Big[\mathcal{N}(s,P_{s-},\Lambda_s,U_s,i)'\mathcal{R}(s,P_{s-},U_s,i)^{-1}\big(P^i_{s-}\sigma^i_s\eta^i_s+\lambda^{\mathbb{G}}_s(P^i_{s-}+U^i_s)F^i_s\gamma^i_s \big)\nonumber\\
&\qquad\qquad\qquad
-(\Lambda^i_s)'\eta^i_s-\lambda^{\mathbb{G}}_sU^i_s\gamma^i_s-\sum_{j=1}^\ell q_{ij}P^j_{s-}(h^j_{s-}-h^i_{s-})\Big]ds\nonumber\\
&\qquad\qquad-\int_{t\wedge\tau}^{T\wedge\tau}(\eta^i_s)'dW_s-\int_{t\wedge\tau}^{T\wedge\tau}\gamma^i_s dM_s, \ \mbox{ for all $i\in\cM$}.
\end{align}
Please note that the coefficients in \eqref{h} are unbounded as so are $\Lambda^i$, $i\in\cM$.

\begin{theorem}\label{hexis}
The system of BSDEs \eqref{h} admits a solution $(h^i,\eta^i,\gamma^i)_{i\in\cM}$ such that
\begin{align*}
(h^i,\eta^i,\gamma^i)\in S^\infty_{\mathbb{G}}(0,T;\mathbb{R})\times L^2_{\mathbb{G}}(0,T;\mathbb{R}^n)\times L^{\infty}_{\mathbb{G}}(0, T;\mathbb{R}), \ \mbox{for all} \ i\in\cM.
\end{align*}
\end{theorem}
\pf
Consider the following system of BSDEs without jumps:
\begin{align}
\label{KBM}
\begin{cases}
dK^{b,i}=\Big[\Big(P^{b,i}\mu^i+\sigma^i\Lambda^{b,i}+\lambda F^i(1-P^{b,i})\Big)'
\Big(P^{b,i}\sigma^i(\sigma^i)'+\lambda G^{a,i}F^i(F^i)'\Big)^{-1}\\
\qquad\qquad\times\Big(K^{b,i}\mu^i+\sigma^i L^{b,i}+\lambda F^i(H^{a,i}-K^{b,i})\Big)
-\lambda(H^{a,i}-K^{b,i})-\sum\limits_{j=1}^{\ell}q_{ij}K^{b,j}
\Big]dt+(L^{b,i})'dW,\\
K^i_{T}=H^{b,i}, \ \mbox{for all} \ i\in\cM.
\end{cases}
\end{align}
By Theorem 3.6 of \cite{HSX2}, the system \eqref{KBM} admits a unique solution $(K^{b,i},L^{b,i})_{i\in\cM}$ such that
\[
(K^{b,i},L^{b,i})\in L^{\infty}_{\mathbb{F}}(0,T;\mathbb{R})\times L^{2}_{\mathbb{F}}(0,T;\mathbb{R}^n), \ \mbox{for all} \ i\in\cM.
\]

Define
\begin{align*}
&K^{i}_t:=K^{b,i}_t\mathbf{1}_{t<\tau}+H^{a,i}_{\tau}\mathbf{1}_{t\geq\tau},\\
&L^{i}_t:=L^{b,i}_t\mathbf{1}_{t\leq\tau},\\
&\zeta^{i}_t:=(H^{a,i}_t-K^{b,i}_t)\mathbf{1}_{t\leq\tau}, \ t\in[0,T], \ i\in\cM.
\end{align*}
By Lemma \ref{link}, $(K^{i},L^{i},\zeta^{i})\in S^{\infty}_{\mathbb{G}}(0,T;\mathbb{R})\times L^2_{\mathbb{G}}(0,T;\mathbb{R}^n)\times L^{\infty}_{\mathbb{G}}(0,T;\mathbb{R})$ is a solution of the following system of BSDEs with jumps:

\begin{align*}
K^i_t&=H^i-\int_{t\wedge\tau}^{T\wedge\tau}\Big[\mathcal{N}(s,P_{s-},\Lambda_s,U_s,i)'
\mathcal{R}(s,P_{s-},U_s,i)^{-1}
\big(K^i_{s-}\mu^i_s+\sigma^i_s L^i_s+\lambda^{\mathbb{G}}_sF^i_s\zeta^i_s\big)\\
&\qquad\qquad\qquad\qquad-\sum_{j=1}^{\ell}q_{ij}K^j_{s-}
\Big]ds-\int_{t\wedge\tau}^{T\wedge\tau}(L^i_s)'dW_s-\int_{t\wedge\tau}^{T\wedge\tau}\zeta^i_s dM_s, \ \mbox{for all} \ i\in\cM.
\end{align*}
Recall that $P^i_{t-}+U^i_t=1$, so we can define $h^i_t=\frac{K^i_t}{P^i_t}, \ \eta^i_t=\frac{L^i_t}{P^i_{t-}}-\frac{h^i_{t-}\Lambda^i_t}{P^i_{t-}}, \ \gamma^i_t=\frac{\zeta^i_t-h^i_{t-}U^i_t}{P^i_{t-}+U^i_t}$, $t\in[0,T]$. Applying It\^{o}'s formula to $\frac{K^i}{P^i}$, it is easy to show that $(h^i,\eta^i,\gamma^i)_{i\in\cM}$ is a solution of \eqref{h}.
\eof

\begin{theorem}
Let $(P^{i},\Lambda^{i},U^{i})_{i\in\cM}$ be the solution of \eqref{Pmv} constructed in Theorem \ref{existence}
and $(h^{i},\eta^{i},\gamma^{i})_{i\in\cM}$ be the solution of \eqref{h} in Theorem \ref{hexis}, then the mean-variance hedging problem \eqref{mvhedge} has an optimal feedback control
\begin{align}\label{optpi}
\pi^*(t,X_{t-},i)=-\mathcal{R}(t,P^i_{t-},U^i_t,i)^{-1}\Big[\mathcal{N}(t,P^i_{t-},\Lambda^i_t,U^i_t,i)(X_{t-}-h_{t-})
-P^i_{t-}\sigma^i_t\eta^i_t-\lambda^{\mathbb{G}}_t\gamma^i_t(P^i_{t-}+U^i_t)F^i_t\Big].
\end{align}
Moreover, the corresponding optimal value is
\begin{align}\label{opthedgevalue}
\min_{\pi\in L^{2}_{\mathbb{H}}(0, T\wedge\tau;\mathbb{R}^m)}\E(X^{\pi}_{T\wedge\tau}-H^{\alpha_{T\wedge\tau}})^2&=P^{i_0}_0(x-h^{i_0}_0)^2+\E\int_0^{T\wedge\tau}\sum_{j=1}^{\ell}q_{\alpha_t j}P^j_{t-}(h^j_{t-}-h^{\alpha_{t-}}_{t-})^2dt\nonumber\\
&\; \quad+\E\int_0^{T\wedge\tau}O(t,\alpha_t)dt,
\end{align}
where
\begin{align*}
O(t,i)&=P^i_{t-}(\eta^i_{t})'\eta^i_{t}+(\gamma^i_{t})^2\lambda^{\mathbb{G}}_{t}
(P^i_{t-}+U^i_{t})-(P^i_{t-}\sigma^i_{t}\eta^i_{t}+\lambda^{\mathbb{G}}_{t}\gamma^i_{t}(P^i_{t-}+U^i_{t})F^i_{t})'\\
& \qquad \times(P^i_{t-}\sigma^i_{t}(\sigma^i_{t})'+\lambda^{\mathbb{G}}_{t}(P^i_{t-}+U^i_{t})F^i_{t}(F^i_{t})')^{-1}
(P^i_{t-}\sigma^i_{t}\eta^i_{t}+\lambda^{\mathbb{G}}_t\gamma^i_{t}(P^i_{t-}+U^i_{t})F^i_{t})\geq 0.
\end{align*}
\end{theorem}
\pf
One can get \eqref{optpi} and\eqref{opthedgevalue} by applying It\^{o}'s formula to $P^{\alpha_t}_t(X_t-h^{\alpha_t}_t)^2$ with some tedious calculation. We shall show $O(t,i)\geq0$. In light of the length of many equations, $``t,t-,i"$ may be suppressed.

From the proof of Theorem \ref{existence}, $P^i_t\geq c>0, \ P^i_{t-}+U^i_t=1$.
If $P\eta'\eta+\gamma^2\lambda^{\mathbb{G}}(P+U)=0$, then $\eta=0$ and $\gamma^2\lambda^{\mathbb{G}}=0$, hence $O(t,i)=0$.

Otherwise, $P\eta'\eta+\gamma^2\lambda^{\mathbb{G}}(P+U)>0$. Then by Lemma \ref{A1}, we have
\begin{align*}
&\qquad\Big[P\eta'\eta+\gamma^2\lambda^{\mathbb{G}}(P+U)
-(P\sigma\eta+\lambda^{\mathbb{G}}\gamma(P+U)F)'
(P\sigma\sigma'+\lambda^{\mathbb{G}}(P+U)FF')^{-1}
(P\sigma\eta+\lambda^{\mathbb{G}}\gamma(P+U)F)\Big]\\
&\qquad\qquad\times\mathrm{det}(P\sigma\sigma'+\lambda^{\mathbb{G}}(P+U)FF')\\
&=(P\eta'\eta+\gamma^2\lambda^{\mathbb{G}}(P+U))\\
&\qquad\times\mathrm{det}\Big(P\sigma\sigma'+\lambda^{\mathbb{G}}(P+U)FF'
-\frac{1}{P\eta'\eta+\gamma^2\lambda^{\mathbb{G}}(P+U)}
(P\sigma\eta+\lambda^{\mathbb{G}}\gamma(P+U)F)(P\sigma\eta+\lambda^{\mathbb{G}}\gamma(P+U)F)'\Big).
\end{align*}
As $P\sigma\sigma'+\lambda^{\mathbb{G}}(P+U)FF'>0$ and $P\eta'\eta+\gamma^2\lambda^{\mathbb{G}}(P+U)>0$,
it suffices to show
\begin{align*}
(P\eta'\eta+\gamma^2\lambda^{\mathbb{G}}(P+U))(P\sigma\sigma'+\lambda^{\mathbb{G}}(P+U)FF')-
(P\sigma\eta+\lambda^{\mathbb{G}}\gamma(P+U)F)(P\sigma\eta+\lambda^{\mathbb{G}}\gamma(P+U)F)'\geq0.
\end{align*}
Actually for any $x\in\mathbb{R}^n$,
\begin{align*}
&\qquad(P\eta'\eta+\gamma^2\lambda^{\mathbb{G}}(P+U))x'(P\sigma\sigma'+\lambda^{\mathbb{G}}(P+U)FF')x\\
&\qquad\qquad-x'(P\sigma\eta+\lambda^{\mathbb{G}}\gamma(P+U)F)(P\sigma\eta+\lambda^{\mathbb{G}}\gamma(P+U)F)'x\\
&=P^2x'\sigma(\eta'\eta I_{n}-\eta\eta')\sigma'x+P\lambda^{\mathbb{G}}(P+U)x'(\gamma\sigma-F\eta')(\gamma\sigma'-\eta F')x\geq 0,
\end{align*}
where we used the fact that $P\lambda^{\mathbb{G}}(P+U)\geq 0$ to get the last inequality.
This completes the proof.
\eof

\begin{remark}
The existence of solutions to \eqref{P}, \eqref{Pmv} and \eqref{h} is proved by construction method (Theorems \ref{existence} and \ref{hexis}), unfortunately we cannot prove the uniqueness so far. We hope to address the uniqueness of solutions to \eqref{P} and \eqref{h} in our future research.
\end{remark}

\begin{appendix}
\section{APPENDIX}\label{appn} 
\begin{lemma}\label{A1}
If $a>0$, $e\in\mathbb{R}^m$, $A\in\mathbb{S}^m$ and $A>0$, then
$$(a-e'A^{-1}e)\;\mathrm{det}(A)=a\;\mathrm{det}(A-\frac{1}{a}ee').$$
\end{lemma}
\pf
By some calculation, we have
\begin{align*}
\begin{bmatrix}
1 & -e'A^{-1}\\
0 & I_m
\end{bmatrix}
\begin{bmatrix}
a & e'\\
e & A
\end{bmatrix}
\begin{bmatrix}
1 & 0\\
-A^{-1}e & I_m
\end{bmatrix}
=
\begin{bmatrix}
a-e'A^{-1}e & 0\\
0 & A
\end{bmatrix},
\end{align*}
and
\begin{align*}
\begin{bmatrix}
1 & 0\\
-\frac{1}{a}e & I_m
\end{bmatrix}
\begin{bmatrix}
a & e'\\
e & A
\end{bmatrix}
\begin{bmatrix}
1 & -\frac{1}{a}e'\\
0 & I_m
\end{bmatrix}
=
\begin{bmatrix}
a & 0\\
0 & A-\frac{1}{a}ee'
\end{bmatrix}.
\end{align*}
Taking determinant on both sides of the two equations, the result follows.
\eof
\end{appendix}
\renewcommand{\refname}

\end{document}